\documentclass[11pt,twoside]{amsart}

\usepackage{geometry}
\usepackage{hyperref}
\usepackage{amsmath,mathrsfs,amssymb,esint,amsthm,bbm,mathtools}

\usepackage[dvips]{graphicx}
\catcode`@=11 \@addtoreset{equation}{section} \catcode`@=12

\def\vbar{\mathchoice{\vrule height6.3ptdepth-.5ptwidth.8pt\kern-.8pt}
   {\vrule height6.3ptdepth-.5ptwidth.8pt\kern-.8pt}
   {\vrule height4.1ptdepth-.35ptwidth.6pt\kern-.6pt}
   {\vrule height3.1ptdepth-.25ptwidth.5pt\kern-.5pt}}

\def\bbc#1#2{{\rm \mkern#2mu\vbar\mkern-#2mu#1}}
\def\bbb#1{{\rm I\mkern-3.5mu #1}}
\def\bba#1#2{{\rm #1\mkern-#2mu{\cal F}udge #1}}
\def\bb#1{{\count4=`#1 \advance\count4by-64 \ifcase\count4\or\bba A{11.5}\or
   \bbb B\or\bbc C{5}\or\bbb D\or\bbb E\or\bbb F \or\bbc G{5}\or\bbb H\or
   \bbb I\or\bbc J{3}\or\bbb K\or\bbb L \or\bbb M\or\bbb N\or\bbc O{5} \or
   \bbb P\or\bbc Q{5}\or\bbb R\or\bba S{8}\or\bba T{10.5}\or\bbc U{5}\or
   \bba V{12}
\or\bba W{16.5}\or\bba X{11}\or\bba Y{11.7}\or\bba Z{7.5}{\cal F}i}}

\newcounter{theorem}
\begin{document}

\bibliographystyle{abbrv}

		\newcommand{\bea}{\begin{eqnarray}}
 		\newcommand{\nvec}{\vec}
 	  \newcommand{\Lform}{\mathscr{L}}
 		\newcommand{\cal}{\mathcal}
    \newcommand{\X} {{\cal X}}
    \newcommand{\ep} {{\varepsilon}}
    \newcommand{\g} {{\bf g}}
    \newcommand{\T}{{\bf T}}
    \newcommand{\F} {{\cal F}}
    \newcommand{\J} {{\cal J}}
    \newcommand{\B} {{\cal B}}
    \newcommand{\Lq} {{\cal L}}
    \newcommand{\cend} {\end{center}}
    \newcommand{\elp}[2]{L^{#1}{(#2)}}
    \newcommand{\jac}{\mathcal{J}}
    \newcommand{\eea}{\end{eqnarray}}
    \newcommand{\rn}{\mathbb{R}^n}
    \newcommand{\fracc}{\displaystyle\frac}
    \newcommand{\intt}{\displaystyle{\int}}
    \newcommand{\ointt}{\displaystyle{\oint}}	
    \newcommand{\disp}{\displaystyle}
    \newcommand{\limm}{\displaystyle\lim}
    \newcommand{\ra}{\rightarrow}
    \newcommand{\sumk} {\sum_{k=1}^{m_j}}
    \newcommand{\sumj} {\sum_{j=1}^n}
    \newcommand{\sumi}{\disp\sum_{i=1}^n}
    \newcommand{\beginc} {\begin{center}}
    \newcommand{\ltwo} {L^2(\Omega)}
    \newcommand{\ltwoj} {{L^2 (\Omega_j)}}
    \newcommand{\cl}[2] {{\cal L}^{#1}_w(#2,{Q})}
    \newcommand{\wa}[2] {W^{1,{#1}}_{w}(#2,{Q})}
    \newcommand{\wb}[2] {W^{1,{#1}}_{w,0}(#2,{Q})}
    \newcommand{\wc}[2] {{\cal W}^{1,{#1}}_{w}(#2,{Q})}
    \newcommand{\we}[2] {{\cal W}^{1,{#1}}_{w,0}(#2,{Q})}
    \newcommand{\normal}[1]{\vec #1}
    \newcommand{\hf} {[f]_{\alpha \; ; \; x_0}}
    \newcommand{\di} {\partial}
    \newcommand{\scq} {\mathcal{Q}}
    \newcommand{\bg} {{\bf g}}
    \newcommand{\diam}{{\rm diam}}
    \newcommand{\ball}[3] {B_{#1}({#2}_{#3})}
    \newcommand{\sumjn} {\sum_{j=1}^N}
    \newcommand{\sumjm} {\sum_{j=1}^M}
    \newcommand{\scl} {\mathcal{L}}
    \newcommand{\rhu} {\rightharpoonup}
    \newcommand{\lhu} {\leftharpoonup}
    \newcommand{\limj} {\disp\lim_{j\rightarrow \infty}}
    \newcommand{\dist}{\textrm{dist}}
    \newcommand{\Q}{{\cal Q}}
    \newcommand{\Po}{{\cal P}}
    
    \newtheorem{remark}[theorem]{{\bf Remark}}
    \newtheorem{defn}[theorem]{{\bf Definition}}
    \newtheorem{thm}[theorem]{\bf{Theorem}}
    \newtheorem{propn}[theorem]{\bf{Proposition}}
    \newtheorem{lemma}[theorem]{\bf{Lemma}}
    \newtheorem{fact}[theorem]{{\bf Fact}}
    \newtheorem{corollary}[theorem]{{\bf Corollary}}

\title [Existence of Weak Solutions]{Existence of Weak Solutions of Linear Subelliptic Dirichlet Problems With Rough Coefficients}
\vspace{-1.0cm}



\maketitle
\begin{center} Scott Rodney\\
Cape Breton University, Nova Scotia Canada
\cend
\vspace{1in}

\begin{center} Abstract
\end{center}

\begin{quote}{\footnotesize This article gives an existence theory for weak solutions of second order non-elliptic linear Dirichlet problems of the form
\bea \nabla'P(x)\nabla u +{\bf HR}u+{\bf S'G}u +Fu &=& f+{\bf T'g} \textrm{ in }\Theta\nonumber\\
u&=&\varphi\textrm{ on }\partial \Theta.\nonumber
\eea
The principal part $\xi'P(x)\xi$ of the above equation is assumed to be comparable to a quadratic form ${\cal Q}(x,\xi) = \xi'Q(x)\xi$ that may vanish for non-zero $\xi\in\mathbb{R}^n$.  This is achieved using techniques of functional analysis applied to the degenerate Sobolev spaces $QH^1(\Theta)=W^{1,2}(\Omega,Q)$ and $QH^1_0(\Theta)=W^{1,2}_0(\Theta,Q)$ as defined in \cite{R1},\cite{SW2},\cite{MRW} and \cite{CRW}.  In \cite{SW1}, with generalizations in \cite{SW2}, the authors give a regularity theory for a subset of the class of equations dealt with here.}
\end{quote}

\section{Notation}

This work will use the basic notation of \cite{SW1} and some of it is recalled now.  Given two subsets $\Theta,\Omega$ of $\mathbb{R}^n$, the notation $\Theta\Subset\Omega$ indicates that $\Theta\subset\Omega$ with $\overline{\Theta}\subset\Omega$. Unless otherwise specified, vectors are always viewed as column vectors.  Given a vector $\xi\in\mathbb{R}^n$, we denote its transpose by $\xi'$.  A vector valued function will always be indicated using bold faced characters except when speaking of the gradient, $\nabla f$, of a function $f$. When two vector valued functions ${\bf H},{\bf G}$ are written as a product ${\bf HG}$, this denotes the standard dot product of ${\bf H}$ with ${\bf G}$.\\

Given a measurable subset $E\subset\mathbb{R}^n$, $|E|$ denotes its Lebesgue measure.  For $N\in\mathbb{N},\;1\leq p<\infty$ $[L^p(E)]^N$ denotes the classical $N$-dimensional Lebesgue space $L^p(E)\times \cdots \times L^p(E)$ with norm 
\bea ||{\bf f}||_{L^p(E)} = ||(f_1,...,f_N)||_{[L^p(E)]^N} = \Big(\disp\sum_{i=1}^N ||f_i||_{L^p(E)}^p\Big)^{1/p}.
\eea
Note that the form of the norm on the right is taken in context with regard to the function ${\bf f}$ inserted on the left.\\

Objects of the form $T(x) = \disp\sum_{i=1}^n t_i(x)\fracc{\partial}{\partial x_i}$ are called vectorfields.  The formal adjoint of $T(x)$ is denoted by $T'(x)$, that is $T'(x)$ is defined by $T'(x)u= -Div( {\bf t}(x)u)$ where ${\bf t}(x)=(t_1(x),...,t_n(x))$.  Given an $N$-tuple ${\bf T} = (T_1(x),...,T_N(x))$ of vectorfields $T_i(x)$ and an $\mathbb{R}^N$-valued function ${\bf g}=(g_1,...,g_N)$, ${\bf Tg}$ denotes the dot product of ${\bf T}$ with ${\bf g}$.  That is,
\bea {\bf Tg} = \disp\sum_{i=1}^N T_i(x)g_i(x) = \disp\sum_{i=1}^N \sum_{j=1}^n t_{i,j}(x)\fracc{\partial}{\partial x_j}g_i(x).
\eea
where $t_{i,j}(x)$ is the $j^{th}$ coefficient of $T_i(x)$.  Similarly, we have
\bea\label{bigsubunit} {\bf T}'{\bf g}=-Div\Big(\disp\sum_{i=1}^N t_{i,1}(x)g_i(x),...,\disp\sum_{i=1}^N t_{i,n}(x)g_i(x)\Big)=-\disp\sum_{j=1}^n\fracc{\partial}{\partial x_j}\disp\sum_{i=1}^Nt_{i,j}(x)g_i(x).
\eea
Lastly, given an $\mathbb{R}^N$-valued function ${\bf H}=(h_1(x),\dots,h_N(x))$ and a real valued function $u(x)$, a term of the form ${\bf HT}u$ denotes 
\bea {\bf H T}u = \disp\sum_{i=1}^N h_i(x)T_i(x)u(x) = \disp\sum_{i=1}^N \sum_{j=1}^n h_i(x)t_{i,j}(x)\fracc{\partial}{\partial x_j}u(x).
\eea
Again, an analogous formula to (\ref{bigsubunit}) holds for terms of the form ${\bf T}'{\bf H}u$.

\section{Introduction}

\cite{SW1} and \cite{SW2} give a regularity theory for weak solutions of second order linear subelliptic equations under the hypothesis of local Poincar\'e/Sobolev type inequalities in the presence of a homogeneous space structure.  This paper gives an existence theory for weak solutions to Dirichlet problems involving second order linear non-elliptic equations of a class containing many of those studied in \cite{SW1}/\cite{SW2}.  Local boundedness estimates and a treatment of Ne\"umann boundary value problems for equations in this class will appear in a sequel.\\

Let $\Omega\subset \mathbb{R}^n$ be open and bounded and $\Theta\subset\Omega$.  Throughout this paper $Q(x)$ denotes a bounded measurable non-negative definite symmetric matrix defined on $\Omega\times \mathbb{R}^n$ and ${\cal Q}(x,\xi)= \xi'Q(x)\xi$ is its related quadratic form.  That is, ${\cal Q}$ is a measurable symmetric quadratic form in $\Omega$ satisfying
\bea \label{A} (i)&& 0\leq {\cal Q}(x,\xi) \textrm{ for all $\xi\in\mathbb{R}^n$ and $a.e.\;x\in\Omega$.  Note that ${\cal Q}(x,\xi)$ may vanish for non-zero }\nonumber\\
&&\xi\in\mathbb{R}^n.\nonumber\\
\label{B} (ii)&& \textrm{There is a $C_0>0$ so that ${\cal Q}(x,\xi) \leq C_0|\xi|^2$ for all $\xi\in\mathbb{R}^n$ and $a.e.\; x\in\Omega$.}\nonumber
\eea
In what follows, this quadratic form ${\cal Q}$ will be referred to as a bounded measurable non-negative definite symmetric quadratic form in $\Omega$.\\

  The focus of this paper is the existence of weak solutions to Dirichlet problems related to linear equations with rough coefficients of the form
\bea\label{pde} Xu=\nabla'P(x)\nabla u +{\bf HR}u+{\bf S'G}u +Fu = f+{\bf T'g} \textrm{ for $x\in \Theta$.}
\eea
Here, $P(x)$ is an $n\times n$ bounded measurable non-negative definite symmetric matrix with related quadratic form $\xi'P(x)\xi$ comparable to ${\cal Q}(x,\xi)$ in $\Theta$.  That is, there are constants $c_1,C_1>0$ so that 
\bea\label{PQcompare}
c_1 \xi'Q(x)\xi \leq \xi'P(x)\xi \leq C_1 \xi'Q(x)\xi
\eea
for all $\xi\in\mathbb{R}^n$ and a.e. $x\in \Theta$.  The coefficient functions ${\bf H,G}$ are, for some $N\in\mathbb{N}$, $\mathbb{R}^N$-valued measurable functions in $\Theta$ and $F$ is a real valued measurable function defined in $\Theta$.  These coefficient functions will belong to Lebesgue spaces whose orders are related to the gain in a given Sobolev inequality, see (\ref{sob}) below.  The data functions $f,{\bf g}$ are, for some $K\in\mathbb{N}$, $L^2(\Theta)$ and $[L^2(\Theta)]^K$ functions respectively.   ${\bf R,S}$ are $N$-tuples and ${\bf T}$ a $K$-tuple of first order vectorfields assumed to be subunit with respect to the quadratic form ${\cal Q}$.  As can be seen in \cite{R2} and \cite{SW1}, a first order vectorfield $V(x) =\sum_{j=1}^n v_j(x)\fracc{\partial}{\partial x_j}$ identified with the vector ${\bf v}(x) = (v_1(x),\dots, v_n(x))$ is said to be subunit with respect to the quadratic form ${\cal Q}$ on $\Theta$ if and only if $v_j$ is defined and measurable in $\Theta$ for each $1\leq j\leq n$ and
\bea\label{subunit}
\Big( {\bf v}(x) \cdot \xi\Big)^2 \leq {\cal Q}(x,\xi)   
\eea
for all $\xi\in \mathbb{R}^n$ and a.e. $x\in \Theta$.  Note that for a function $w\in Lip(\Theta)$, $V(x)w(x)={\bf v}(x)\cdot \nabla w(x)$ for almost every $x\in \Theta$ and (\ref{subunit}) gives
\bea \Big( {\bf v}(x) \cdot \nabla w(x)\Big)^2 \leq {\cal Q}(x,\nabla w(x)) \textrm{ for almost every $x\in \Theta$}.\nonumber
\eea

As in the elliptic case presented in \cite{GT}, a negativity condition for the lower order terms ${\bf G,S}$ and $F$ of $X$ will be required.  That is, it will be assumed that
\bea\label{negativity} 
\intt_{\Theta} (Fw + {\bf GS}w)dx \leq 0
\eea
for all $w\in Lip_0(\Theta)$ satisfying $w(x)\geq 0$ in $\Theta$. 
\begin{remark} As will be seen below in the proof of Theorem \ref{th1}, condition (\ref{negativity}) is the key property that allows the application of the Fredholm Alternative enabling one to conclude existence of weak solutions.  This can also be seen in the elliptic case.  For example, setting ${\bf G}={\bf H}=\vec{0},\;{\bf g}=\nvec{0}$, $F=c$ for a fixed constant $c$ and $P(x)=Id$, equation (\ref{pde}) becomes the elliptic equation
\bea \Delta u +cu = f.
\eea
Here, the negativity condition (\ref{negativity}) becomes $c\leq 0$ which is sufficient for the existence of weak solutions to equations of this type, see \cite{GT}.  Lastly, condition (\ref{negativity}) may seem to differ from that of \cite{GT} by a negative sign, see \cite{GT} condition (8.8), but they are equivalent.  This is due to the usage of the formal adjoint ${\bf S}'$ of the vectorfield ${\bf S}$ in (\ref{pde}).  This term appears as $-S'$ in \cite{GT}.
\end{remark}

In current literature, the word "subelliptic", when applied to equations of the form (\ref{pde}), refers to second order equations possessing certain regularity properties, (see \cite{SW1}).  These equations are, up to condition (\ref{negativity}), a special case of the equations one can solve via the methods presented in this work.  This motivates the following definition naming operators that belong to what is referred to here as "the subelliptic class".  This definition will be used in the statements of the main theorems greatly reducing the length of their statements.

\begin{defn} \label{pdeclass} Let $\Theta\subset\Omega$.  A second order operator $X$ of the form 
\bea\label{Xclass} X = \nabla'P(x)\nabla + {\bf HR} + {\bf S'G} + F
\eea
is said to be of the subelliptic class related to $(\Omega,{\cal Q},\Theta)$ if and only if
\bea &&(i)\;\;\; P(x) \textrm{ is a bounded measurable non-negative definite symmetric matrix defined in $\Theta$}\nonumber\\
&&\;\;\;\;\;\;\; \textrm{satisfying (\ref{PQcompare}),}\nonumber\\
&&(ii) \;\;{\bf R,S} \textrm{ are, for some $N\in \mathbb{N}$, $N$-tuples of first order vectorfields subunit with respect to}\nonumber\\
&&\;\;\;\;\;\;\; \textrm{${\cal Q}$ in $\Theta$,}\nonumber\\
&&(iii)\; {\bf H,G} \textrm{ are measurable $\mathbb{R}^N$-valued functions defined in $\Theta$, $F$ is a real valued measura-}\nonumber\\
&&\;\;\;\;\;\;\;\textrm{ble function defined in $\Theta$ and}\nonumber\\
&&(iv) \;\textrm{${\bf S},{\bf G},F$ satisfy the negativity condition (\ref{negativity}).}\nonumber
\eea
\end{defn}

The major hypotheses of the results contained in this article are local Poincar\'e and Sobolev inequalities adapted to the quadratic form ${\cal Q}$.  These inequalities will be assumed to hold on quasimetric balls given by a quasimetric $\mu(x,y)$ defined in $\Omega$ and upper semicontinuous in the second variable.  The quasimetric ball of radius $r>0$ centred at $x\in \Omega$ is given by
\bea \label{dball} B_r(x) = \{y\in\Omega\; :\; \mu(x,y)<r\}.
\eea
For the purposes of this work, it will be assumed that the pair $(\Omega,\mu)$ is a homogeneous space.  As in \cite{SW1}, a pair $(\Omega,\mu)$ is a homogeneous space if $\mu$ is as above and the collection of quasimetric balls $\{B_r(y)\}_{r>0;y\in\Omega}$ satisfies a doubling condition with respect to Lebesgue measure.  That is, there are constants $c_2>1,C_2>0$ so that
\bea\label{double} |B_{c_2r}(y)| \leq C_2 |B_r(y)|
\eea
for all $y\in \Omega$ and $r>0$.  The reader is referred to section 2.2 of \cite{SW1} for a development of properties of homogeneous spaces.  In the case where ${\cal Q}(x,\cdot)$ is continuous in $\Omega$, one can picture the quasimetric $\mu$ as the Carnot-Carath\'eodory control metric for familiarity, however the homogeneous space structure is the essential property needed for this work.  We assume that there are constants $C_3,\;C_4,\;\delta>0$ and $\sigma >1$ so that the following inequalities hold on these balls.  Note that for $x\in\Omega$, the quantity dist$(x,\partial\Omega)$ denotes the Euclidean distance from $x$ to the boundary of $\Omega$.\\

\noindent{\bf The Poincar\'e Inequality:}\\
Given $x\in\Omega$ and $0<r<\delta\textrm{dist}(x,\partial\Omega)$ we have that
\bea\label{pc} 
\Big(\fint_{B_r(x)} |w(y)-w_{B_r(x)}|^2dy\Big)^{1/2} &\leq& C_3r\Big(\fint_{B_r(x)}{\cal Q}(y,\nabla w)dy\Big)^{1/2}
\eea
for all $w\in Lip_{\cal Q}(B_r(x))=\{\varphi\in Lip(B_r(x))\; :\; ||\varphi||_{L^2(B_r(x))}+||{\cal Q}(x,\nabla \varphi)||_{L^1(B_r(x))}^{1/2}<\infty\}$ where $w_{B_r(x)}=\fint_{B_r(x)} w(y)dy$ is the average of $w$ over $B_r(x)$.\\

\noindent {\bf The Sobolev Inequality:}\\
Given $x\in \Omega$ and $0<r<\delta\textrm{dist}(x,\partial\Omega)$ we have that 
\bea\label{sob} 
\Big(\fint_{B_r(x)} |w(y)|^{2\sigma}dy\Big)^{1/2\sigma} &\leq& C_4\Big[\Big(\fint_{B_r(x)} |w(y)|^2dy\Big)^{1/2} + r\Big(\fint_{B_r(x)} {\cal Q}(y,\nabla w)dy\Big)^{1/2}\Big]
\eea
for all $w\in Lip_0(B_r(x))$.\\

These local inequalities give rise to an inequality of the utmost importance for this work and is presented in the following proposition.\\

\begin{propn}\label{prelim01} Let $\Omega$ be a bounded open subset of $\mathbb{R}^n$ and suppose that the pair $(\Omega,\mu)$ defines a homogeneous space in the sense of \cite{SW1}.  Let $\Theta\Subset \Omega$ be open and suppose that the local Poincar\'e and Sobolev inequalities, (\ref{pc}) and (\ref{sob}) respectively, hold. Then, the global Sobolev inequality 
\bea\label{globsob} \Big(\intt_\Theta |w(x)|^{2\sigma}dx\Big)^{1/2\sigma} \leq C_5\Big(\intt_\Theta {\cal Q}(x,\nabla w(x))dx\Big)^{1/2}
\eea
holds for all $w\in Lip_0(\Theta)$ where $C_5>0$ is independent of $w$.
\end{propn}
\noindent This proposition is proved in \cite{R2} using the existence of weak solutions to a parametrized family of Dirichlet problems related to operators in the subelliptic class satisfying ${\bf H}={\bf G}=\vec{0}$.  This development is quite lengthy and the reader is referred to \cite{R2} for a complete presentation.  

\begin{remark}  The proof of Proposition \ref{prelim01} in \cite{R2} is given with $\mu$ identified with the Carnot-Carath\'eodory control metric.  It should be noted that nowhere in the argument is it needed that $\mu$ be so specific.  Rather, the argument uses that $\mu$ is a quasimetric for which the pair $(\Omega,\mu)$ is a homogeneous space.
\end{remark}

The global Sobolev inequality (\ref{globsob}) given by Proposition \ref{prelim01} is the key to solving Dirichlet problems related to operators in the subelliptic class.  As the local Sobolev and Poincar\'e inequalities are also needed for Proposition (\ref{compact}) below, the theorems presented in this paper will contain the hypotheses of Proposition \ref{prelim01} in their statements. \\

\begin{remark} 
In the elliptic case, where $Q(x,\xi)=|\xi|^2$, the classical Sobolev inequality has the form (\ref{globsob}), for $n\geq 3$, where $\sigma=\fracc{n}{n-2}>1$ and $C_5=\fracc{2(n-1)}{\sqrt{n}(n-2)}$, see \cite{GT}.
\end{remark}

Before presenting the main theorems, a definition of Sobolev space and weak solution is required.  The development of the degenerate Sobolev spaces presented here is quite brief.  For a complete discussion see \cite{SW2} and also \cite{CRW} for weighted versions.     

\begin{defn} Let $\Theta\subset\Omega$ and ${\cal Q}$ a measurable non-negative definite symmetric quadratic form defined in $\Omega$.\\
{\bf (i)}  The degenerate Sobolev space $QH^1(\Theta)=W^{1,2}(\Theta,Q)$ is defined as the completion of $Lip_{\cal Q}(\Theta)=\{\varphi\in Lip(\Theta) : ||\varphi||_{L^2(\Theta)}+\intt_\Theta {\cal Q}(x,\nabla\varphi)dx<\infty\}$ with respect to the norm
\bea\label{norm} ||w||_{QH^1(\Theta)} &=& ||w||_{L^2(\Theta)} + ||\nabla w||_{{\cal L}^2(\Theta,Q)}.
\eea
{\bf (ii)}  The degenerate Sobolev space $QH^1_0(\Theta)=W^{1,2}_0(\Theta,Q)$ is defined as the completion of $Lip_{{\cal Q},0}(\Theta)=\{\varphi\in Lip_0(\Theta) : \intt_\Theta {\cal Q}(x,\nabla\varphi)dx<\infty\}$ with respect to the norm (\ref{norm}).\\
\noindent {\bf (iii)} Let $\alpha\in \mathbb{R}$.  The degenerate Sobolev space $QH^1_\alpha(\Theta)$ is defined as the completion of $Lip_{{\cal Q},\alpha}(\Theta)=\{\varphi\in Lip(\Theta) : ||\varphi||_{L^2(\Theta)}+\intt_\Theta {\cal Q}(x,\nabla\varphi)dx<\infty \textrm{ and } \overline{supp(\varphi-\alpha)}\subset\Theta\}$ with respect to the norm (\ref{norm}).

\end{defn}

\begin{remark}
\begin{enumerate}
\item The space ${\cal L}^2(\Theta,Q)$ is defined in \cite{SW2} as the collection of equivalence classes of measurable vector valued functions ${\bf f}(x)=(f_1(x),\dots,f_n(x))$ for which 
\bea \label{scriptLnorm} || {\bf f} ||_{{\cal L}^2(\Theta,Q)} &=& \Big(\intt_\Theta {\cal Q}(x,{\bf f}(x))dx\Big)^{1/2} < \infty.
\eea
Two vector valued functions ${\bf f,g}$ are regarded as equivalent if and only if $||{\bf f} - {\bf g}||_{{\cal L}^2(\Theta,Q)} = 0$.  It is shown in \cite{SW2} that ${\cal L}^2(\Theta,Q)$ is a Hilbert space with respect to the inner product
\bea \label{innerp} <{\bf f},{\bf g}>_{\cal L} &=& \intt_\Theta ({\bf f}(x))'Q(x){\bf g}(x)dx.
\eea
\item In view of functional analysis methods; using that $\Theta$ is bounded, H\"older's inequality together with (\ref{globsob}) gives that the norm (\ref{scriptLnorm}) is equivalent to the norm (\ref{norm}) on $QH^1_0(\Theta)$. 
\end{enumerate}
\end{remark}

As defined, $QH^1_0(\Theta)$ consists of equivalence classes of Cauchy sequences of Lipschitz functions with compact support in $\Theta$.  For a sequence $\{f_j\}\in QH^1_0(\Theta)$, its equivalence class is denoted by $[\{f_j\}]$.  Since $\{f_j\}$ is Cauchy in $QH^1_0(\Theta)$, it is also Cauchy in $L^2(\Theta)$ and the sequence $\{\nabla f_j\}$ is Cauchy in ${\cal L}^2(\Theta,Q)$.  As $L^2(\Theta)$ and ${\cal L}^2(\Theta,Q)$ are Banach spaces, to the equivalence class $[\{f_j\}]$ one assigns a unique pair $(f,{\bf g})\in L^2(\Theta)\times {\cal L}^2(\Theta,Q)$.  This gives a Banach space isomorphism ${\cal J}:QH^1_0(\Theta) \ra E\subset L^2(\Theta)\times {\cal L}^2(\Theta,Q)$ defined by 
\bea \label{isomorph}{\cal J}([\{f_j\}]) = (f,{\bf g}).
\eea
The image $E$ of $QH^1_0(\Theta)$ under ${\cal J}$ is denoted in \cite{SW2} by ${\cal W}^{1,2}_0(\Theta,Q)$ and is a closed subset of $L^2(\Theta)\times {\cal L}^2(\Theta,Q)$.  Since these spaces are isomorphic, no distinction will be made between them.  For simplicity, the vector function ${\bf g}$ in (\ref{isomorph}) will be referred to as $\nabla f$ even though ${\bf g}$ may not be uniquely determined by $f$; see the example found in \cite{FKS} and/or \cite{R1} for a development of this idea.  Further abuse of notation will be made when referring to an element of $QH^1_0(\Theta)$.  Indeed, the notation $u\in QH^1_0(\Omega)$ refers to a specific pair $(u,{\bf g})\in L^2(\Theta)\times {\cal L}^2(\Theta,Q)$ (written as $(u,\nabla u)$) and/or the equivalence class of Cauchy sequences of Lipschitz functions written, again with abuse of notation, as $[u]$.  When working with a specific sequence of $[u]$, it will be referred to as a sequence representing $u$.  Lastly, as can be seen in \cite{R1} and \cite{SW2}, $QH^1_0(\Theta)$ is a Hilbert space with respect to the inner product
\bea\label{inner2} ((u,\nabla u),(v,\nabla v)) = \intt_\Theta uv + \intt_\Theta (\nabla u)'Q(x)\nabla v.
\eea

\begin{remark}\label{extendsob}
\begin{enumerate}
\item  If inequality (\ref{globsob}) holds it is evident that (\ref{globsob}) extends by density to hold for all $w\in QH^1_0(\Theta)$ with the same constant $C_5$.  That is, for an element $w=(w,{\bf g})\in QH^1_0(\Theta)$, the extended Sobolev inequality is
\bea \Big(\intt_\Theta |w(x)|^{2\sigma}dx\Big)^{1/2\sigma} \leq C_5\Big(\intt_\Theta {\cal Q}(x,{\bf g}(x))dx\Big)^{1/2}.
\eea
This inequality will be used in the form (\ref{globsob}) with the abuse of notation ${\bf g}=\nabla w$.
\item  In the case where ${\cal Q}(x,\cdot)$ is bounded (i.e. the case where $Q(x)$ is a bounded non-negative definite symmetric measurable matrix defined in $\Omega$), $Lip_{{\cal Q},0}(\Theta) = Lip_0(\Theta)$.  Thus, in the elliptic setting (when ${\cal Q}(x,\xi)=|\xi|^2$), $QH^1_0(\Theta) = H^1_0(\Theta)$.
\end{enumerate}
\end{remark}

In order to define a notion of weak solution of equation (\ref{pde}), we turn to the divergence formula for motivation.  For a suitable set $\Theta\subset\Omega$, setting {\bf n} as the unit outward normal to $\partial \Theta$, the divergence theorem says that 
\bea\label{divform} \;\;\;\;\;\intt_{\partial\Theta}  \Big(vP(x)\nabla u\Big)\cdot {\bf n}=\intt_\Theta \textrm{div}\Big(vP(x)\nabla u\Big)  = \intt_\Theta (\nabla v)'P(x)\nabla u + \intt_\Theta v\Big((\nabla )'P(x)\nabla\Big) u
\eea
for $u\in C^2(\overline{\Theta}),\;v\in C^1(\overline{\Theta})$. Hence, for $v$ with compact support in $\Theta$ this simplifies to 

\bea\label{divform2}
\intt_\Theta v\Big((\nabla )'P(x)\nabla\Big) u = -\intt_\Theta (\nabla v)'P(x)\nabla u .
\eea
Consider now the bi-linear form $\Lform:QH^1(\Theta)\times QH^1_0(\Theta)\ra \mathbb{R}$ associated to equation (\ref{pde}), motivated by (\ref{divform2}), defined by 
\bea\label{bilinear} \Lform(u,v) = \intt_\Theta (\nabla v)'P(x)\nabla u - \intt_\Theta v{\bf HR}u - \intt_\Theta u{\bf GS}v - \intt_\Theta Fuv.
\eea

\begin{defn}\label{weaksol1}
\begin{enumerate}\item An element $u\in QH^1_0(\Theta)$ is called a weak solution of the homogeneous Dirichlet problem
\bea \label{dp} Xw &=& f+{\bf T'g} \textrm{ in }\Theta\nonumber\\
w&=&0\textrm{ on }\partial\Theta
\eea
if and only if the equality 
\bea \label{weaksol}
\Lform(u,v) &=& -\intt_\Theta fv - \intt_\Theta {\bf gT}v
\eea
holds for every $v\in QH^1_0(\Theta)$.
\item A function $u\in QH^1(\Theta)$ is called a weak solution of the Dirichlet problem
\bea \label{dpnh}Xw &=& f+{\bf T'g} \textrm{ in }\Theta\nonumber\\
w&=&\varphi\textrm{ on }\partial\Theta
\eea 
if and only if $u$ satisfies (\ref{weaksol}), $\varphi\in QH^1(\Theta)$ and $u-\varphi=(u-\varphi,\nabla u - \nabla \varphi)\in QH^1_0(\Theta)$.
\end{enumerate}
\end{defn}

We now state the first result of this paper.  The reader is reminded here that ${\cal Q}(x,\cdot)$ is a bounded measurable non-negative definite symmetric quadratic form defined in $\Omega\times \mathbb{R}^n$.

\begin{thm} \label{th1} Suppose that the pair $(\Omega,\mu)$ defines a homogeneous space in the sense of \cite{SW1} and that the local Poincar\'e and Sobolev inequalities (\ref{pc}) and (\ref{sob}) hold.  Fix $q>2\sigma'$, where $\sigma'$ is the dual of the exponent $\sigma$ appearing in (\ref{sob}), and let $\Theta\Subset\Omega$ be open.  Fix $N,K\in\mathbb{N}$ and let $X$ be an operator of the subelliptic class related to $(\Omega,{\cal Q},\Theta)$ with coefficients ${\bf H},{\bf G}\in [L^q(\Theta)]^N$ and $F\in L^{q/2}(\Theta)$. If $f\in L^2(\Theta)$, ${\bf g}\in [L^2(\Theta)]^K$ and ${\bf T}$ is a $K$-tuple of first order vectorfields subunit with respect to ${\cal Q}$ in $\Theta$ then the homogeneous Dirichlet problem (\ref{dp}) is uniquely solvable in $QH^1_0(\Theta)$.
\end{thm}

As a corollary to this, we obtain an existence result for the general Dirichlet problem
\bea\label{dp2} Xu &=& f+{\bf T'g} \textrm{ in }\Theta\nonumber\\
u&=&\varphi\textrm{ on }\partial\Theta.
\eea

\begin{thm}\label{thm3}Suppose that the pair $(\Omega,\mu)$ defines a homogeneous space in the sense of \cite{SW1}.  Assume that the Poincar\'e inequality (\ref{pc}) and the Sobolev inequality (\ref{sob}) hold. Let $\Theta\Subset\Omega$ be an open set, fix $q,r,s\geq 1$ so that $q>2\sigma'$ and $q\geq \max\{2r,2s\}$.  Fix $N,K\in \mathbb{N}$ and let $X$ be an operator of the subelliptic class related to $(\Omega,{\cal Q},\Theta)$ with coefficients ${\bf H}, {\bf G}\in [L^{q}(\Theta)]^N$ and $F\in L^{q}(\Theta)$. Assume that the data $f \in L^2(\Theta), {\bf g}\in [L^2(\Theta)]^K$ and that ${\bf T}$ is a $K$-tuple of first order vectorfields subunit with respect to ${\cal Q}$ in $\Theta$.  Further, suppose that the boundary value $\varphi\in QH^1(\Theta)$ with $\varphi\in L^{2r'}(\Theta)$ and ${\cal Q}(x,\nabla\varphi)\in L^{s'}(\Theta)$. Then, the Dirichlet problem (\ref{dp2}) is uniquely solvable in $QH^1(\Theta)$.  
\end{thm}

The reader may feel at this point that the integrability conditions imposed on $(\varphi,\nabla \varphi)$ seem artificial.  However, Proposition \ref{prelim01} offers some answers.  Indeed, assume that $\Theta \subset {\cal V}\subset \Omega$ are open sets with $\overline{{\cal V}}\subset\Omega$ and that $\varphi\in QH^1_0({\cal V})$.  If $(\Omega,\mu)$ is a homogeneous space in the sense of \cite{SW1} and the local Sobolev and Poincar\'e inequalities hold, Proposition \ref{prelim01} gives that the boundary value $\varphi\in L^{2\sigma}({\cal V})$ and hence, as $\Theta\subset{\cal V}$, $\varphi \in L^{2\sigma}(\Theta)$.  The next theorem incorporates this with Theorem \ref{thm3}.  Since the proof of this next theorem is identical to the proof of the previous upon setting $r=\sigma'$, it is omitted.

\begin{thm}\label{th4} Suppose that the pair $(\Omega,\mu)$ defines a homogeneous space in the sense of \cite{SW1}.  Assume that the Poincar\'e inequality (\ref{pc}) and the Sobolev inequality (\ref{sob}) hold.  Let $\Theta \subset {\cal V}\subset \Omega$ be open sets with $\overline{\cal V}\subset \Omega$.  Fix $q,s\geq 1$ so that $q>2\sigma'$ and $q\geq 2s$.  Let $N,K\in\mathbb{N}$ and $X$ be an operator of the subelliptic class related to $(\Omega,{\cal Q},\Theta)$ with coefficients ${\bf H},{\bf G}\in [L^q(\Theta)]^N$ and $F\in L^q(\Theta)$. Assume that the data $f \in L^2(\Theta), {\bf g}\in [L^2(\Theta)]^K$ and that ${\bf T}$ is a $K$-tuple of first order vectorfields subunit with respect to ${\cal Q}$ in $\Theta$.   Further, suppose that the boundary value $\varphi\in QH^1_0({\cal V})$ with ${\cal Q}(x,\nabla \varphi)\in L^{s'}(\Theta)$. Then, the Dirichlet problem (\ref{dp2}) is uniquely solvable in $QH^1(\Theta)$.
\end{thm}

\section{Preliminaries}

In order to shorten proofs, the following collection of preliminary propositions and lemmas is given.  The first result to be discussed is a compact "embedding" result for degenerate Sobolev spaces.  The proof of this theorem can be found in \cite{R2} and with weighted versions in \cite{CRW}.

\begin{propn}\label{compact} Let $\Omega$ be a bounded open subset of $\mathbb{R}^n$ and $(\Omega,\mu)$ a homogeneous space in the sense of \cite{SW1}.  Let $\Theta$ be an open subset of $\Omega$ with $\overline{\Theta}\subset \Omega$.  Let $Q$ be a bounded measurable non-negative semidefinite symmetric matrix in $\Omega$ for which the local Poincar\'e inequality (\ref{pc}) and local Sobolev inequality (\ref{sob}) hold.  Then the projection map $P:QH^1_0(\Omega)\ra  L^2(\Theta)$ defined by 
\bea\label{proj} P((u,\nabla u)) = u
\eea
is a compact mapping.
\end{propn}
\begin{remark} The projection map in Proposition \ref{compact} is not referred to as an embedding as it fails to be 1-1 in general.  See \cite{R1} and \cite{SW2} for more details.
\end{remark}

Before focusing on the bilinear form $\Lform$, it is essential that terms of the form $Tu$ and $Tuv$ for $u,v\in QH^1_0(\Theta)$ are defined when $T$ is a vectorfield subunit to ${\cal Q}$.  This is achieved via the action of a subunit vectorfield on $Lip_Q(\Theta)$ functions.  

\begin{lemma}\label{subunitstuff} Suppose that $\Theta$ is a bounded subset of $\mathbb{R}^n$.  Let $(u,\nabla u),(v,\nabla v)\in QH^1(\Theta)$ with representative sequences $\{u_j\},\{v_j\}$.  Then, if $T = {\bf t}(x)\cdot\nabla $ is subunit with respect to ${\cal Q}$ in $\Theta$ we have that\\

\noindent 1.  there is a unique $L^2(\Theta)$ function associated to $(u,\nabla u)$ that we denote by $Tu={\bf t}(x)\cdot \nabla u$.  Further, if $u\in Lip_{\cal Q}(\Theta)$ then $Tu={\bf t}(x)\cdot\nabla u$ where $\nabla u$ is the gradient of $u$ defined almost everywhere in $\Theta$.\\

\noindent 2.  If $(u,\nabla u),(v,\nabla v)\in QH^1_0(\Theta)$ and the global Sobolev inequality (\ref{globsob}) holds, there is a unique $L^\frac{2\sigma}{\sigma+1}(\Theta)$ function $Tuv$ associated to $(u,\nabla u)$ and $(v,\nabla v)$.  Moreover, $Tuv = uTv+vTu$ as $L^\frac{2\sigma}{\sigma+1}(\Theta)$ functions where $Tu,Tv$ are defined as in 1.
\end{lemma}

\begin{remark} Note that the notation $Tu$ is an abuse since $T$ is acting on $u=(u,{\bf g})\in QH^1(\Theta)$.  Formally, $Tu = T(u,{\bf g})={\bf t}(x)\cdot{\bf g}$ is accurate but we will employ $Tu$ for notational convenience.
\end{remark}

\noindent{\bf Proof:}  For each $j$, set $w_j(x) = Tu_j(x)$.  Then, $w_j$ is a measurable function since $u_j\in Lip(\Theta)$ and $T$ has measurable coefficients.  Also, $\{w_j\}$ is Cauchy in $L^2(\Theta)$.  Indeed,
\bea ||w_i-w_j||_{L^2(\Theta)}&=& ||T(u_i-u_j)||_{L^2(\Theta)}\nonumber\\
&\leq& ||{\cal Q}(x,\nabla u_i - \nabla u_j)||_{L^1(\Theta)}^{1/2}\textrm{ (by subuniticity)}\nonumber\\
&\leq& ||u_i-u_j||_{QH^1(\Theta)}.
\eea
Thus, $\{w_j\}$ is Cauchy in $L^2(\Theta)$ as $\{u_j\}$ is Cauchy in $QH^1(\Theta)$.  Therefore, a unique function $w\in L^2(\Theta)$ is associated to $(u,\nabla u)$ and we write $w= Tu$.  This gives 1.\\

Next, suppose that $(u,\nabla u),(v,\nabla v)\in QH^1_0(\Theta)$ and that the global Sobolev inequality (\ref{globsob}) holds.  Notice that for each $j\in \mathbb{N}$, $Tu_jv_j = u_jTv_j + v_jTu_j$ a.e. in $\Theta$.  Thus, setting $s = \fracc{2\sigma}{\sigma +1}$ we have by adding and subtracting $u_iTv_j + v_iTu_j$,
\bea ||Tu_iv_i - Tu_jv_j||_{L^s(\Theta)}&\leq & ||u_iT(v_i-v_j)||_{L^s(\Theta)} + ||v_iT(u_i-u_j)||_{L^s(\Theta)}\nonumber\\
&& +||(u_i-u_j)T(v_j)||_{L^s(\Theta)}
+ ||(v_i-v_j)T(u_j)||_{L^s(\Theta)}\nonumber\\
&=& I+II+III+IV.
\eea
Terms I and II of the above are estimated in the same way as is the case for terms III and IV.  Therefore, we focus only on terms I and III.  Using H\"older's inequality with exponent $\sigma +1$ and conjugate $\frac{\sigma+1}{\sigma}$, we obtain
\bea \label{Tuv1}I^s &\leq& ||u_i||_{L^{2\sigma}(\Theta)}^s ||  T(v_i-v_j)||_{L^2(\Theta)}^s\nonumber\\
&\leq&C_5^s||u_i||_{QH^1(\Theta)}^s ||v_i-v_j||_{QH^1(\Theta)}^s
\eea
by the Sobolev inequality (\ref{globsob}) and using that $T$ is subunit with respect to ${\cal Q}$ in $\Theta$.  Using the same argument, we have 
\bea \label{Tuv2}III^s&\leq&C_5^s||u_i-u_j||_{QH^1(\Theta)}^s ||v_j||_{QH^1(\Theta)}^s.
\eea
Therefore, noting similar estimates for II and IV, $\{Tu_iv_i\}$ is Cauchy in $L^s(\Theta)$. Thus, to the pairs $(u,\nabla u),(v,\nabla v)$ we associate the unique function $Tuv = \disp\lim_{j\ra\infty} Tu_jv_j$ where the limit is taken in $L^s(\Theta)$.  Lastly, to see that $Tuv = uTv+vTu$ as $L^s(\Theta)$ functions, recall that $Tu_jv_j = u_jTv_j + v_jTu_j$ a.e.  Further, the sequences $\{u_iTv_i\}$ and $\{v_iTu_i\}$ both converge in $L^s(\Theta)$ to $uTv$ and $vTu$ respectively by the same application of H\"older's inequality above.  Thus, $Tuv = uTv+vTu$ in $L^s(\Theta)$.  $\Box$\\

\begin{remark}\label{subunitremark} Lemma \ref{subunitstuff} indicates that integrals of the form $\intt_\Theta HTuv dx$ are finite provided $H\in L^q(\Theta)$ where $q\geq 2\sigma'$ since $s' = \Big(\fracc{2\sigma}{\sigma +1}\Big)^{'} = 2\sigma'$.  Further, we may write $\intt_\Theta HTuvdx = \intt_\Theta HvTudx + \intt_\Theta HuTvdx$ for $u,v\in QH^1_0(\Theta)$.
\end{remark} 

To simplify computations we give the following lemma demonstrating a pointwise estimate for terms of the form $|{\bf T}u|$ where ${\bf T}=(T_1,...,T_N)$ with each $T_i$ a vectorfield subunit with respect to ${\cal Q}$ in $\Theta$.

\begin{lemma}\label{suprop} Fix $N\in\mathbb{N}$ and ${\bf T}(x) = (T_1(x),...,T_N(x))$ where each $T_i(x) = \sum_{j=1}^n t_{i,j}(x) \fracc{\partial}{\partial x_j}$ identified with the vector ${\bf t}(x)=(t_{i,1},...,t_{i,n})$ satisfies 
\bea \Big({\bf t}(x)\cdot \xi \Big)^2 \leq Q(x,\xi)
\eea
for every $\xi\in \mathbb{R}^n$ and a.e. $x\in \Theta$.  Then, for any $(u,\nabla u)\in QH^1_0(\Theta)$ and a.e. $x\in\Theta$ we have the pointwise estimate
\bea \label{vectorsubestimate}|{\bf T}(x)u|&\leq&\sqrt{N}\sqrt{{\cal Q}(x,\nabla u)}.
\eea
\end{lemma}

\noindent {\bf Proof:} Fix $(u,\nabla u)\in QH^1(\Theta)$. The proof of this result is the following calculation using subuniticity of ${\bf T}$.

\bea |{\bf T}(x)u| = \Big(\sum_{i=1}^N \Big(T_i(x) u\Big)^2\Big)^{1/2}&=&  \Big(\sum_{i=1}^N \Big({\bf t}(x)\cdot \nabla u\Big)^2\Big)^{1/2}\nonumber\\
&\leq& \Big(\sum_{i=1}^N {\cal Q}(x,\nabla u)\Big)^{1/2}= \sqrt{N}\sqrt{{\cal Q}(x,\nabla u)} \textrm{ by (\ref{subunit}).}\;\; \Box\nonumber
\eea

In order to apply functional analysis techniques to the bilinear form $\Lform$ acting on the Hilbert space $QH^1_0(\Omega)$ we require boundedness and coercivity estimates for $\Lform$.  These are contained in the next two lemmas.

\begin{lemma}\label{prelim02} Suppose that $\Theta\subset\Omega$ and that the global Sobolev inequality (\ref{globsob}) holds.  Fix $N\in\mathbb{N}$ and suppose that ${\bf R}=(R_1,...,R_N),{\bf S}=(S_1,...,S_N)$ are $N$-tuples of first order vector fields subunit with respect to ${\cal Q}$ in $\Theta$, ${\bf H},{\bf G}\in [L^{2\sigma'}(\Theta)]^N$ and that $F\in L^{\sigma'}(\Theta)$ where $\sigma'$ is the dual exponent to $\sigma>1$ in (\ref{globsob}).  Then, the bilinear form $\Lform$ satisfies
\bea \label{bounded} |\Lform (u,v)| \leq C||u||_{QH^1(\Theta)}||v||_{QH^1(\Theta)}
\eea
for all $u,v\in QH^1_0(\Theta)$ where $C>0$ is independent of both $u$ and $v$.
\end{lemma}

\noindent{\bf Proof:} Let $u,v\in QH^1_0(\Theta)$.  Applying Lemma \ref{suprop} to ${\bf R}$ and ${\bf S}$ we have
\bea |\Lform (u,v)|&\leq& \intt_\Theta |(\nabla u)'P(x)(\nabla v)|dx   + \intt_\Theta |v{\bf HR}u|dx + \intt_\Theta |u{\bf GS}v|dx + \intt_\Theta |Fuv|dx\nonumber\\
&\leq& C_1^2\intt_\Theta |{\cal Q}(x,\nabla u)|^{1/2}|{\cal Q}(x,\nabla v)|^{1/2}dx + \sqrt{N}\intt_\Theta |v{\bf H}||{\cal Q}(x,\nabla u)|^{1/2}dx\nonumber\\
&& + \sqrt{N}\intt_\Theta |u{\bf G}||{\cal Q}(x,\nabla v)|^{1/2}dx +\intt_\Theta |F|^{1/2}|u||F|^{1/2}|v|dx\nonumber\\
&\leq& C_1^2||u||_{QH^1(\Theta)}||v||_{QH^1(\Theta)} + \sqrt{N}||u||_{QH^1(\Theta)}\Big(\intt_\Theta |v{\bf H}|^2dx\Big)^{1/2} \nonumber\\
&&+ \sqrt{N}||v||_{QH^1(\Theta)}\Big(\intt_\Theta |u{\bf G}|^2dx\Big)^{1/2} + \Big(\intt_\Theta |F||u|^2dx\Big)^{1/2}\Big(\intt_\Theta |F||v|^2dx\Big)^{1/2}\nonumber.
\eea
This last line is obtained via an application of the Cauchy Schwarz inequality.  Applying H\"older's inequality with exponents $\sigma, \sigma'$ we obtain
\bea |\Lform (u,v)|&\leq&C_1^2||u||_{QH^1(\Theta)}||v||_{QH^1(\Theta)} + \sqrt{N}||{\bf H}||_{L^{2\sigma'}(\Theta)}||u||_{QH^1(\Theta)}||v||_{L^{2\sigma}(\Theta)}\nonumber\\
&& + \sqrt{N}||{\bf G}||_{L^{2\sigma'}(\Theta)}||v||_{QH^1(\Theta)}||u||_{L^{2\sigma}(\Theta)} + ||F||_{L^{\sigma'}(\Theta)}||u||_{L^{2\sigma}(\Theta)}||v||_{L^{2\sigma}(\Theta)}.
\eea

Lastly, we apply the Sobolev inequality (\ref{globsob}) (noting Remark \ref{extendsob}-1) to the $L^{2\sigma}(\Theta)$ norms appearing above giving 

\bea |\Lform (u,v)|&\leq& C||u||_{QH^1(\Theta)}||v||_{QH^1(\Theta)}
\eea
where $C=(C_1^2 +\sqrt{N}C_5(||{\bf H}||_{L^{2\sigma'}(\Theta)} + ||{\bf G}||_{L^{2\sigma'}(\Theta)}) + C_5^2 ||F||_{L^{\sigma'}(\Theta)})$.  $\Box$\\

The next lemma yields the "almost coercivity" of the bilinear form $\Lform$.  This result indicates that the operator $X$ is only a compact operator away from having an associated coercive bilinear form.  In fact, as will be seen in the remark following the lemma, if the lower order terms of $X$ satisfy ${\bf H}={\bf G}={\vec 0} , F=0$, then $\Lform$ is coercive and existence of weak solutions to the homogeneous Dirichlet problem (\ref{dp}) is concluded via the Lax-Milgram theorem.  

\begin{lemma}\label{coercive} Suppose that $\Theta\subset\mathbb{R}^n$ is open and that the global Sobolev inequality (\ref{globsob}) holds.  Fix $N\in\mathbb{N},q>2\sigma'$ and suppose that ${\bf R}=(R_1,...,R_N),{\bf S}=(S_1,...,S_N)$ are $N$-tuples of first order vectorfields subunit with respect to ${\cal Q}$ in $\Theta$, that ${\bf H},{\bf G}\in [L^{q}(\Theta)]^N$ and that $F\in L^{q/2}(\Theta)$.  Then, the bilinear form $\Lform$ satisfies the almost coercive estimate
\bea\label{almost} |\Lform(u,u)| \geq \fracc{c_1}{4}||u||_{QH^1(\Theta)}^2 - C||u||_{L^2(\Theta)}^2
\eea
for all $u\in QH^1_0(\Omega)$ where $c_1>0$ is as in (\ref{PQcompare}) and $C>0$ is independent of $u$.
\end{lemma}

\noindent {\bf Proof:} Fix $u\in QH^1_0(\Theta)$.  By the triangle inequality we have that $\Lform$ satisfies
\bea \label{I-II}|\Lform(u,u)| &\geq& c_1\intt_\Theta {\cal Q}(x,\nabla u)dx - |\intt_\Theta u({\bf HR}+{\bf GS})udx| - |\intt_\Theta Fu^2dx|\nonumber\\
&=& c_1\intt_\Theta {\cal Q}(x,\nabla u)dx - I - II.
\eea
We estimate terms $I$ and $II$ separately beginning with $I$.  By Lemma (\ref{suprop}) we have that
\bea\label{I} |I|&\leq& \sqrt{N}\intt_\Theta |u({\bf H}+{\bf G})||{\cal Q}(x,\nabla u)|^{1/2}dx\nonumber\\
&\leq&  \sqrt{N}\Big(\intt_\Theta u^2|{\bf H}+{\bf G}|^2dx\Big)^{1/2}\Big(\intt_\Theta {\cal Q}(x,\nabla u)dx\Big)^{1/2}.\nonumber\\
&\leq & \sqrt{N}||{\bf H}+{\bf G}||_{L^{q}(\Theta)}||u^2||_{L^{\frac{q}{q-2}}(\Theta)}^{1/2}\Big(\intt_\Theta {\cal Q}(x,\nabla u)dx\Big)^{1/2}
\eea
where the last line is obtained using an application of H\"older's inequality with $\fracc{q}{2}>\sigma'$ and $\Big(\fracc{q}{2}\Big)^{'}=\fracc{q}{q-2}$.  Since $1<\Big(\fracc{q}{2}\Big)^{'}<\sigma$, we have by the interpolation inequality of \cite{GT}, chapter 7, that
\bea \label{interpI} ||u^2||_{L^\frac{q}{q-2}(\Theta)} &\leq& C^2_{\epsilon,q,\sigma}||u^2||_{L^1(\Theta)} + \epsilon^2 ||u^2||_{L^\sigma(\Theta)}\nonumber\\
&=& C^2_{\epsilon,q,\sigma}||u||_{L^2(\Theta)}^2 + \epsilon^2 ||u||_{L^{2\sigma}(\Theta)}^2
\eea
where $C_{\epsilon,q,\sigma}=\epsilon^{-\frac{\sigma}{q(\sigma-1)-2\sigma}}$ and $\epsilon>0$ is to be chosen in a moment.  Note that $\fracc{\sigma}{q(\sigma-1)-2\sigma}>0$ as $q>2\sigma'$.  Set $M=||{\bf H}||_{L^q(\Theta)}+||{\bf G}||_{L^q(\Theta)}+ ||F||_{L^{q/2}(\Theta)}$.  Using that $(a+b)^{1/2} \leq a^{1/2}+b^{1/2}$ for $a,b\geq 0$, inserting (\ref{interpI}) into (\ref{I}) with an application of the Sobolev inequality (\ref{globsob}) gives
\bea \label{ac-2} |I|&\leq& \sqrt{N} C_{\epsilon,q,\sigma}M||u||_{L^2(\Theta)}\Big[\intt_\Theta {\cal Q}(x,\nabla u)dx\Big]^{1/2} + \sqrt{N}\epsilon C_5M\intt_\Theta {\cal Q}(x,\nabla u)dx\nonumber\\
&\leq& \fracc{2}{c_1}NC_{\epsilon,q,\sigma}^2M^2||u||^2_{L^2(\Theta)} + \Big(\fracc{c_1}{2} +\sqrt{N}\epsilon C_5M\Big)\intt_\Theta {\cal Q}(x,\nabla u)dx.
\eea
Choosing $\epsilon=\fracc{c_1}{8\sqrt{N}C_5M}$ we have that 
\bea\label{Ilast} |I|&\leq& C_I||u||^2_{L^2(\Theta)} + \fracc{5c_1}{8}\intt_\Theta {\cal Q}(x,\nabla u)dx
\eea
where $C_I=\fracc{2}{c_1}NC_{\epsilon,q,\sigma}^2M^2$.\\

Turning our attention to term $II$ of (\ref{I-II}), H\"older's inequality with exponents $\frac{q}{2}, \frac{q}{q-2}$ gives the following using the same argument as in (\ref{interpI}).
\bea\label{II-1} |II|&\leq& ||F||_{L^{q/2}(\Theta)}||u^2||_{L^\frac{q}{q-2}(\Theta)}\nonumber\\
&\leq& C^2_{\eta,q,\sigma}M||u^2||_{L^1(\Theta)} + \eta^2 M||u^2||_{L^\sigma(\Theta)}\nonumber\\
&=& C^2_{\eta,q,\sigma}M||u||^2_{L^2(\Theta)} + \eta^2 M||u||^2_{L^{2\sigma}(\Theta)}.
\eea
Setting $\eta = \sqrt{\frac{c_1}{8C_5^2M}}$, the global Sobolev inequality (\ref{globsob}) gives
\bea \label{IIlast} |II|&\leq& C_{II}||u||_{L^2(\Theta)}^2 + \fracc{c_1}{8}\intt_\Theta {\cal Q}(x,\nabla u)dx
\eea 
where $C_{II}=C_{\eta,q,\sigma}^2M^2$.\\

Inserting (\ref{Ilast}) and (\ref{IIlast}) into (\ref{I-II}) yields
\bea \label{last} |\Lform (u,u)| \geq \fracc{c_1}{4}\intt_\Theta {\cal Q}(x,\nabla u)dx - C_{III} ||u||_{L^2(\Theta)}^2
\eea
where $C_{III} = C_I+C_{II}$.  Adding and subtracting $\fracc{c_1}{4}||u||_{L^2(\Theta)}^2$ on the right gives (\ref{almost}) with $C=C_{III}+\fracc{c_1}{4}$.  $\Box$\\

\begin{remark} Since $\Theta$ is a bounded subset of $\mathbb{R}^n$, if ${\bf H}={\bf G}=\vec{0}, \;F=0$ we can show that ${\Lform}$ is a coercive bilinear form via an application of the Sobolev inequality (\ref{globsob}).  Indeed, for $(w,\nabla w)\in QH^1_0(\Theta)$ we have that 
\bea ||w||_{QH^1(\Theta)}^2 &=& ||w||_{L^2(\Theta)}^2 + ||\nabla w||_{{\cal L}^2(\Theta)}^2\nonumber\\
&\leq& \Big(1+|\Theta|^{\frac{1}{\sigma'}}C_5^2\Big)||\nabla w||_{{\cal L}^2(\Theta)}^2\nonumber\\
&\leq&\Big(\fracc{1+|\Theta|^{\frac{1}{\sigma'}}C_5^2}{c_1}\Big){\Lform}(w,w).
\eea
Hence, dividing by $\fracc{1}{c_1}(1+|\Theta|^{\frac{1}{\sigma'}}C_5^2)$ gives that $\Lform$ is coercive on $QH^1_0(\Theta)$.  This leads immediately to the following existence result via the Lax-Milgram theorem.  

\begin{thm}  Suppose that the pair $(\Omega,\mu)$ defines a homogeneous space in the sense of \cite{SW1} and that the local Poincar\'e and Sobolev inequalities (\ref{pc}) and (\ref{sob}) hold.  Fix $q>2\sigma'$, where $\sigma'$ is the dual of the exponent $\sigma$ appearing in (\ref{sob}), and let $\Theta\Subset\Omega$ be open.  Fix $K\in\mathbb{N}$ and let $X$ be an operator of the subelliptic class related to $(\Omega,{\cal Q},\Theta)$ with coefficients ${\bf H},{\bf G}=\vec{0}$ and $F=0$. Then, with data $f\in L^2(\Theta),\; {\bf g}\in [L^2(\Theta)]^K$ and ${\bf T}=(T_1,...,T_K)$ with each $T_i$ a subunit vector field with respect to ${\cal Q}$ in $\Theta$, the homogeneous Dirichlet problem (\ref{dp}) is uniquely solvable in $QH^1_0(\Theta)$.
\end{thm}
\end{remark}

We now turn to the Homogeneous Dirichlet problem with data zero and show that if it admits a weak solution then that solution must be zero a.e.  In order to prove such a result, a definition of $w_+$ for a generic element $w\in QH^1_0(\Theta)$ is needed.  This requires a chain rule for $QH^1_0(\Theta)$ functions.  The next lemma is an adaptation to $QH^1_0(\Theta)$ of lemma 20 found in \cite{SW2} and the the reader is referred there for its proof.  

\begin{lemma}\label{SWlemma}(Sawyer $\&$ Wheeden 2008)  Suppose that $(u,\nabla u)\in QH^1(\Theta)$ where $\Theta$ is a bounded subset of $\mathbb{R}^n$.  \\

\noindent 1.  If $\phi\in Lip_{{\cal Q}}(\Theta)$ then the pair $(\phi u,\phi \nabla u + u\nabla \phi)\in QH^1(\Theta)$.\\

\noindent 2.  If $f\in C^1(\mathbb{R})$ with $f'\in L^\infty(\mathbb{R})$ then $(f\circ u, (f'\circ u)\nabla u)\in QH^1(\Theta)$. Also, if $u(x)\geq a >b$ and $f\in C^1((b,\infty))$ with $f'\in L^\infty((b,\infty))$ then $(f\circ u, (f'\circ u)\nabla u)\in QH^1(\Theta)$.\\

\noindent 3.  Both $\Big(u_+,{\cal X}_{\{u>0\}}\nabla u\Big)$ and $\Big(u_-,{\cal X}_{\{u<0\}}\nabla u\Big)$ are elements of $QH^1(\Theta)$.  Moreover, if $(u,\nabla u)\in QH^1_\alpha(\Theta)$ with $\alpha\geq 0$ then $\Big(u_+,{\cal X}_{\{u>0\}}\nabla u\Big)\in QH^1_\alpha(\Theta)$ and $\Big(u_-,{\cal X}_{\{u<0\}}\nabla u\Big)\in QH^1_0(\Theta)$ - the opposite holds if $\alpha<0$.\\

\noindent ** For a set $E$, ${\cal X}_E$ denotes the characteristic function of $E$.
\end{lemma}

\begin{remark} Although not relevant in the current scheme, item 2 and 3 above do hold for unbounded $Q$; See \cite{SW2} Lemma 20.
\end{remark}

\begin{propn}\label{zeroprob}  Suppose that $\Theta\subset\Omega$ is open and that the global Sobolev inequality (\ref{globsob}) holds.  Fix $q>2\sigma'$, $N\in\mathbb{N}$ and let $X$ be of the subelliptic class related to $(\Omega,{\cal Q},\Theta)$.  Assume that the coefficients ${\bf H},{\bf G}\in [L^{q}(\Theta)]^N$ and that $F\in L^{q/2}(\Theta)$.  If $u\in QH^1_0(\Theta)$ is a weak solution of the Dirichlet problem 
\bea\label{dpzero}Xw &=& 0 \textrm{ in }\Theta\nonumber\\
w&=&0 \textrm{ on }\partial \Theta
\eea
then $u=0$ in $QH^1_0(\Theta)$, that is, $u=0$ a.e. in $\Theta$ and $\nabla u =0$ in ${\cal L}^2(\Theta,Q)$.
\end{propn}

\noindent{\bf Proof:} We begin by showing that $u=0$ a.e. in $\Theta$.  Let $(u,\nabla u)\in QH^1_0(\Theta)$ be a weak solution of the Dirichlet problem (\ref{dpzero}) such that $u\neq 0$ and let $\{u_j\}$ be a sequence representing $(u,\nabla u)$.  Note that since $u\neq 0$, the Sobolev inequality (\ref{globsob}) gives that $supp({\cal Q}(x,\nabla u))$ has positive measure.  Choose $k>0$ so that $0<k<\disp\sup_\Theta u$ (If no such $k$ exists, see the paragraph following (\ref{last})).  Then, the pair $(u-k,\nabla u)\in QH^1_{-k}(\Theta)$.  Indeed, set $w_j = u_j - k$ for each $j\in\mathbb{N}$.  Then $w_j\in Lip_{Q,-k}(\Theta)$. Clearly, $\nabla w_j=\nabla u_j$ a.e. and $||(u-k) - w_j||_{L^2(\Theta)} = ||u-u_j||_{L^2(\Theta)}$.  Hence, $\nabla w_j \ra \nabla u$ in ${\cal L}^2(\Theta)$ and $w_j\ra u-k$ in $L^2(\Theta)$.  We conclude that the pair $(u-k,\nabla u)\in QH^1_{-k}(\Theta)$.  By part 3. of Lemma \ref{SWlemma}, we have that $v=((u-k)_+,{\cal X}_{u>k}\nabla u)\in QH^1_0(\Theta)$ is a valid test function for $u$.  That is, the following equality holds:
\bea \label{zero1} \Lform(u,v) = \intt_\Theta (\nabla u)'P(x)\nabla vdx - \intt_\Theta v{\bf HR}u dx - \intt_\Theta u{\bf GS}vdx - \intt_\Theta Fuvdx=0.
\eea
Using that ${\bf GS}(uv) =  u{\bf GS}v + v{\bf GS}u$ (see Remark \ref{subunitremark}) and that $\nabla u = \nabla v$ on the support of $\nabla v$, we have by the negativity condition (\ref{negativity}) that
\bea \intt_\Theta (\nabla v)'P(x)\nabla v dx \leq \intt_\Theta \Big[v{\bf HR}v - v{\bf GS}v\Big]dx \leq \intt_\Theta |v|\Big[|{\bf H}||{\bf R}v| + |{\bf G}||{\bf S}v|\Big]dx
\eea
Using that the vectorfields $R_i$ and $S_i$ are subunit with respect to ${\cal Q}$ in $\Theta$ and that $c_1{\cal Q}(x,\xi)\leq \xi'P(x)\xi$ for each $\xi\in\mathbb{R}^n$ and a.e. $x\in \Theta$, we have with $\Gamma = supp \Big({\cal Q}(x,\nabla v)\Big)$ that 
\bea\label{firstgamma} \Big(\intt_\Theta {\cal Q}(x,\nabla v)dx\Big)^{1/2} &\leq& \fracc{\sqrt{N}}{c_1}||v(|{\bf H}|+|{\bf G}|)||_{L^2(\Gamma)}.
\eea
We estimate the right hand side of (\ref{firstgamma}) using H\"older's inequality with exponents $\frac{q}{2}$ and $\frac{q}{q-2}$.  Indeed,
\bea ||v(|{\bf H}|+|{\bf G}|)||_{L^2(\Gamma)} \leq ||v||_{L^\frac{2q}{q-2}(\Gamma)}\Big(||{\bf H}||_{L^q(\Theta)}+||{\bf G}||_{L^q(\Theta)}\Big) = C({\bf H},{\bf G})||v||_{L^\frac{2q}{q-2}(\Gamma)}.
\eea
Since $q>2\sigma'$ we use H\"older's inequality once again with exponent $p= \frac{\sigma(q-2)}{q}>1 $ and dual $p'=\frac{\sigma(q-2)}{\sigma(q-2)-q}$ giving
\bea ||v||_{L^\frac{2q}{q-2}(\Gamma)} \leq |\Gamma|^{\frac{1}{2\sigma'} - \frac{1}{q}}||v||_{L^{2\sigma}(\Gamma)}\leq |\Gamma|^{\frac{1}{2\sigma'} - \frac{1}{q}}||v||_{L^{2\sigma}(\Theta)}.
\eea
Finally, applying the Sobolev inequality (\ref{globsob}) we obtain the estimate
\bea \label{HGestimate} ||v||_{L^\frac{2q}{q-2}(\Gamma)}\leq C_5|\Gamma|^{\frac{1}{2\sigma'} - \frac{1}{q}}\Big(\intt_\Theta {\cal Q}(x,\nabla v)dx\Big)^{1/2}.
\eea
Combining (\ref{HGestimate}) with (\ref{firstgamma}) gives
\bea \Big(\intt_\Theta {\cal Q}(x,\nabla v)dx\Big)^{1/2} \leq \fracc{\sqrt{N}C_5C({\bf H},{\bf G})}{c_1} |\Gamma|^{\frac{1}{2\sigma'}-\frac{1}{q}}\Big(\intt_\Theta {\cal Q}(x,\nabla v)dx\Big)^{1/2}.
\eea
Dividing by $\Big(\intt_\Theta {\cal Q}(x,\nabla v)dx\Big)^{1/2}$, noting that $\frac{1}{2\sigma'} - \frac{1}{q}>0$ as $q>2\sigma'$ by hypothesis, yields
\bea \label{last} |\Gamma| \geq \Big( \fracc{c_1}{\sqrt{N}C_5C({\bf H},{\bf G})}\Big)^{\frac{2q\sigma'}{q-2\sigma'}}>0.
\eea
This estimate is independent of $k$ and so sending $k\ra \sup u$ gives us that $v=0$ a.e. while the $L^1(\Theta)$ function $(\nabla v)'Q\nabla v$ has support with positive measure.  This contradicts inequality (\ref{firstgamma}). \\

Repeating the above argument with $v=((u+k)_-,{\cal X}_{u<-k}\nabla u)$, where $k>0$ is chosen so that $0<k<-\inf u$, gives that $u=0$ a.e. in $\Theta$.  \\

Suppose again that $u\in QH^1_0(\Theta)$ is a weak solution of (\ref{dpzero}).  Then, $u=(0,{\bf y})$ for some ${\bf y}\in {\cal L}^2(\Theta,Q)$.  Further, for any $v\in QH^1_0(\Theta)$ we have
\bea \Lform(u,v) &=& \intt_\Theta (\nabla u)'P(x)\nabla vdx - \intt_\Theta v{\bf HR}u dx - \intt_\Theta u{\bf GS}vdx - \intt_\Theta Fuvdx\nonumber\\
&=& \intt_\Theta ({\bf y})'P(x)\nabla vdx- \intt_\Theta v{\bf HR}u dx\nonumber\\
&=&0\nonumber
\eea
and so 
\bea \label{grad0}\intt_\Theta ({\bf y}(x))'P(x)\nabla v(x)dx &=& \intt_\Theta v(x){\bf H}(x){\bf R}(x)u dx\nonumber\\
& =&\intt_\Theta v(x){\bf H}(x)\cdot \Big(\sumj r_{1,j}(x)y_j(x),...,\sumj r_{N,j}(x)y_j(x)\Big) dx
\eea
where $y_j$ is the $j^{th}$ component of ${\bf y}$, $R = (R_1,...,R_N)$ and $r_{i,j}$ is the $j^{th}$ coefficient of $R_i$.  Since (\ref{grad0}) holds for every $v\in QH^1_0(\Theta)$ it must also hold for $v=u=(0,{\bf y})$.  That is,
\bea \label{grad0-1}  \intt_\Theta ({\bf y}(x))'P(x){\bf y}(x) dx = 0.
\eea
Thus, using (\ref{PQcompare}), it must be the case that ${\bf y}=0$ in ${\cal L}^2(\Theta,Q)$ giving that $u=0$ in $QH^1_0(\Theta)$.  $\Box$\\   

\section{Proof of Theorem \ref{th1}}

We begin by remarking that Proposition \ref{prelim01} applies and we have a global Sobolev inequality of the form ({\ref{globsob}) holding on $\Theta$.  Recall that we are interested in finding a weak solution of the Dirichlet problem
\bea\label{dp-proof} Xw&=& f+{\bf T}'{\bf g}\textrm{ in } \Theta,\\
w&=&0\textrm{ on }\partial \Theta.\nonumber
\eea
That is, we wish to find an element $u=(u,\nabla u)\in QH^1_0(\Theta)$ so that 
\bea \label{L-1}
{\Lform}(u,v) = -\intt_\Theta fv - \intt_\Theta {\bf gT}v
\eea
holds for every $v\in QH^1_0(\Theta)$. \\

For $p>0$, consider the bilinear form ${\Lform}_p:QH^1_0(\Theta)\times QH^1_0(\Theta)\ra \mathbb{R}$ defined by 
\bea\label{p-form} {\Lform}_p(w,v) = {\Lform}(w,v) + p\intt_\Theta wv.
\eea
Fix $p>C$ where $C$ is as in Lemma \ref{coercive}.  Lemmas \ref{prelim02} and \ref{coercive} then give that ${\Lform}_p$ is both bounded and coercive on the Hilbert space $QH^1_0(\Theta)$.  Thus, the Lax-Milgram theorem gives that for any $\F\in \Big(QH^1_0(\Theta)\Big)^*$ there is a unique $y=(y,\nabla y)\in QH^1_0(\Theta)$ for which 
\bea\label{L-2} {\Lform_p}(y,v) = \F(v)
\eea
for every $v=(v,\nabla v)\in QH^1_0(\Theta)$. That is, the equation 
\bea\label{L-3} X_p w = \F
\eea
is uniquely solvable in $QH^1_0(\Theta)$ where $X_p = X - p$.  Hence, $X_p^{-1}$ is a continuous 1-1 mapping of $\Big(QH^1_0(\Theta)\Big)^*$ onto $QH^1_0(\Theta)$. \\

In order to avoid confusion in the next development, an element $(w,\nabla w)\in QH^1_0(\Theta)$ will be denoted by ${\bf w}$.  For each ${\bf w}=(w,\nabla w)\in QH^1_0(\Theta)$ define a map $J:QH^1_0(\Theta)\ra \Big(QH^1_0(\Theta)\Big)^*$ by setting
\bea\label{L-6} J{\bf w}((v,\nabla v))=\intt_\Theta wv.
\eea
Then $J$ is a compact mapping.  Indeed, we may write $J = I_1P$ where $P:QH^1_0(\Theta)\ra L^2(\Theta)$ is the mapping of Proposition \ref{compact} defined by $P({\bf w})=w$ and $I_1:L^2(\Theta)\ra \Big(QH^1_0(\Theta)\Big)^*$ is defined by $I_1w({\bf v})=I_1w(v,\nabla v)=\intt_\Theta wv$.  As $I_1$ is continuous and $P$ compact, the composition $I_1P=J$ is compact.  \\

To conclude the proof, notice that the original equation $X{\bf w} = f+{\bf T}'{\bf g}$ is equivalent to 
\bea\label{L-5} X_p{\bf w} + pJ{\bf w} = \F
\eea
where $\F$ is the continuous linear functional defined by $\F((v,\nabla v)) =  -\intt_\Theta fv - \intt_\Theta {\bf gT}v$.  The continuity of ${\cal F}$ is due to the fact that ${\bf T}$ is a $K$-tuple of first order vectorfields subunit with respect to ${\cal Q}$ in $\Theta$.  Applying $X_p^{-1}$ to (\ref{L-5}) we obtain
\bea \label{L-7} {\bf w} + pX_p^{-1}J{\bf  w} = X_p^{-1}\F.
\eea
The mapping $-pX_p^{-1}J$ is compact by the compactness of $J$ and the continuity of $X_p^{-1}$.  Therefore, the existence of a unique ${\bf u}=(u,\nabla u)\in QH^1_0(\Theta)$ satisfying (\ref{L-7}) is guaranteed by Proposition 3.23 and the Fredholm Alternative.  $\Box$

\section{Proof of Theorem \ref{thm3}}

 Assume the hypotheses and consider the homogeneous Dirichlet problem
\bea\label{newdp} Xu &=& f+\T'\g - \nabla'P(x)\nabla \varphi - {\bf HR}\varphi - {\bf S}'{\bf G}\varphi - F\varphi \textrm{ in }\Theta,\nonumber\\
u&=&0 \textrm{ on } \partial\Theta.
\eea
The goal is to show that this equation satisfies the hypotheses of Theorem \ref{th1}.  Set $f_1 = {\bf H}{\bf R}\varphi,\; {\bf f_2} = {\bf G}\varphi \textrm{ and } f_3=F\varphi$.  Then, since $\varphi\in L^{2r'}(\Theta)$ and ${\cal Q}(x,\nabla \varphi)\in L^{s'}(\Theta)$ we have that
\bea\label{l2} (i)&&f_1\in L^2(\Theta), \\
\label{l3}(ii)&& {\bf f_2}\in [L^2(\Theta)]^N, \textrm { and }\\
\label{l4}(iii)&&f_3\in L^2(\Theta).
\eea 
To verify (\ref{l2}) we use H\"older's inequality and the subuniticity of the components of {\bf R}.  Indeed,
\bea\label{l5} ||f_1||_{L^2(\Theta)}^2 = \intt_\Theta |{\bf HR}\varphi|^2 &\leq& N\intt_\Theta |{\bf H}|^2{\cal Q}(x,\nabla \varphi)\nonumber\\
&\leq& N|||{\bf H}|^2||_{L^s(\Theta)}||{\cal Q}(x,\nabla \varphi)||_{L^{s'}(\Theta)}\nonumber\\
&=& N||{\bf H}||_{L^{2s}(\Theta)}^2||{\cal Q}(x,\nabla \varphi)||_{L^{s'}(\Theta)}  <\infty
\eea
since ${\cal Q}(x,\nabla\varphi)\in L^{s'}(\Theta)$, $q\geq 2s$ and $\Theta$ is bounded.  Verification of (\ref{l3}) is achieved using that $\varphi\in L^{2r'}(\Theta)$ and $q\geq 2r$.  Indeed,
\bea\label{l6} ||{\bf f_2}||_{L^2(\Theta)}^2 &=& \intt_\Theta |{\bf G}\varphi|^2\nonumber\\
&\leq& |||{\bf G}|^2||_{L^r(\Theta)}||\varphi^2||_{L^{r'}(\Theta)}\nonumber\\
&=& ||{\bf G}||_{L^{2r}(\Theta)}^2||\varphi||_{L^{2r'}(\Theta)}^2<\infty
\eea
since $q\geq 2r$ and $\Theta$ is bounded.  (\ref{l4}) is verified using the same argument as for (\ref{l3}) replacing ${\bf G}$ with $F$.\\

Before moving forward we concern ourselves with the term $-\nabla'P(x)\nabla \varphi$.  Since $P(x)$ is non-negative definite and symmetric for each $x\in \Theta$, it has a unique symmetric non-negative definite square root $\sqrt{P(x)}$ defined in $\Theta$, see appendix 2 of \cite{MRW}.  With this, we may write $\xi'P(x)\xi = (\sqrt{P(x)}\xi)'\sqrt{P(x)}\xi$ for every $\xi\in\mathbb{R}^n$.  Thus, the term of the bilinear form associated to $-\nabla'P(x)\nabla \varphi$ can be written as 
\bea \label{Q} -\intt_\Theta (\nabla\varphi)'P(x)\nabla vdx &=& -\intt_\Theta (\sqrt{P(x)}\nabla \varphi )' \sqrt{P(x)}\nabla vdx\\
&=&-\intt_\Theta (\sqrt{C_1}\sqrt{P(x)}\nabla \varphi )'(\fracc{1}{\sqrt{C_1}} \sqrt{P(x)}\nabla v)dx.\nonumber
\eea  
This integrand is viewed as the dot product of the $[L^2(\Theta)]^n$ function ${\bf f_4}=-\sqrt{C_1}\sqrt{P(x)}\nabla \varphi$ with an $n$-tuple of subunit vectorfields ${\bf V}(x)=(V_1(x),...,V_n(x))$ applied to $(v,\nabla v)$.  Here, $V_i(x)=\fracc{1}{\sqrt{C_1}}{\bf p}_i(x)\cdot \nabla$ where ${\bf p}_i(x)$ is the $i^{th}$ row vector of $\sqrt{P(x)}$.  Thus we write $-\nabla'P(x)\nabla \varphi = {\bf V}'{\bf f_4}$ where ${\bf V}'$ is the adjoint of ${\bf V}$.  To see that each $V_i(x)$ is subunit with respect to ${\cal Q}$ in $\Theta$ we perform the following calculation.  Let ${\bf \xi}\in\mathbb{R}^n$, $\xi_j$ denote its $j^{th}$ coordinate and $p_{i,j}$ denote the $j^{th}$ coordinate of ${\bf p}_i$.  Then, we have for each $1\leq i\leq n$ that

\bea\label{subforrootQ} {\bf \xi}'Q(x){\bf \xi} &\geq& \fracc{1}{C_1}{\bf \xi}'P(x){\bf \xi}\nonumber\\
&=&\fracc{1}{C_1}(\sqrt{P(x)}{\bf \xi})'(\sqrt{P(x)}{\bf \xi})\nonumber\\
&=& \fracc{1}{C_1}\Big(\sumj p_{1,j}(x)\xi_j,...,\sumj p_{n,j}(x)\xi_j\Big)\cdot \Big(\sumj p_{1,j}(x)\xi_j,...,\sumj p_{n,j}(x)\xi_j\Big)\nonumber\\
&=& \fracc{1}{C_1}\disp\sum_{i=1}^n\Big(\sumj p_{i,j}(x)\xi_j\Big)^2\nonumber\\
&\geq& \fracc{1}{C_1}\Big(\sumj p_{i,j}(x)\xi_j\Big)^2\nonumber\\
&=&\Big(\fracc{1}{\sqrt{C_1}}{\bf p}_i(x)\cdot {\bf{ \xi}}\Big)^2.
\eea
Thus, $V_i(x)$ is subunit with respect to ${\cal Q}$ in $\Theta$ as claimed.

With these identifications we rewrite (\ref{newdp}) in the equivalent form
\bea\label{newdp1} X u &=& \overline{f} + \overline{{\bf T}}'{\overline {\bf g}}\textrm{ in $\Theta$},\nonumber\\
u&=&0\textrm{ on $\partial\theta$}
\eea
where $\overline{f} = f - f_1 - f_3\in L^2(\Theta)$, $\overline{\bf g}=({\bf g},-{\bf f_2},{\bf f_4})\in [L^2(\Theta)]^{K+N+n}$ with 
\bea {\overline{\bf T}}=(T_1,...,T_K,S_1,...,S_N,V_1,...,V_n)\nonumber
\eea
a $K+N+n$ tuple of vectorfields subunit to ${\cal Q}$ in $\Theta$. \\

We now apply Theorem \ref{th1} to obtain a unique weak solution $u\in QH^1_0(\Theta)$ of (\ref{newdp1}).  Setting $w=u+\varphi\in QH^1(\Theta)$ gives a weak solution of (\ref{dp2}) by the linearity of $X$.  To see that our weak solution $w$ is unique, notice that if $(y,\nabla y)\in QH^1(\Theta)$ is a weak solution of $(\ref{dp2})$ then by definition $(y-\varphi,\nabla y - \nabla\varphi)\in QH^1_0(\Theta)$ is a weak solution of (\ref{newdp}).  Thus, $(y-\varphi,\nabla y - \nabla\varphi)=(u,\nabla u)$ giving that $y=w$ in $QH^1(\Theta)$.  $\Box$


\end{document}